\documentclass[11pt,leqno,oneside]{amsart}
\usepackage[width=6.5in,height=8.8in]{geometry}
\usepackage[mathscr]{eucal}
\usepackage{amssymb, amsmath,array, amscd, bbm, comment, charter}
\usepackage{enumerate}
\usepackage{graphicx}
\usepackage{url}
\usepackage[colorlinks,plainpages,backref]{hyperref}

\newtheorem{thm}{Theorem}[section]

\newtheorem{ex}[thm]{Example}
\theoremstyle{definition}
\newtheorem{rem}[thm]{Remark}

\begin{document}
\title{Upper level sets of Lelong numbers on Hirzebruch surfaces}

\author{Al\.i Ula\c{s} \"Ozg\"ur K\.i\c{s}\.isel}
\address{Department  of Mathematics,  Middle East  Technical University,  06800 Ankara,  Turkey}
\email{akisisel@metu.edu.tr}

\author{Ozcan Yazici}
\address{Department  of Mathematics,  Middle East  Technical University,  06800 Ankara,  Turkey} 
\email{oyazici@metu.edu.tr}
\date{\today}

\begin{abstract} Let $\mathbb F_a$ denote the Hirzebruch surfaces and $\mathcal{T}_{\alpha,\alpha^{\prime}}(\mathbb{F}_{a})$ denotes the set of positive, closed $(1,1)$-currents on $\mathbb{F}_{a}$ whose cohomology class is $\alpha F+\alpha^{\prime} H$ where $F$ and $H$ generates the Picard group of $\mathbb F_a$. $E^+_{\beta}(T)$ denotes the upper level sets of Lelong numbers $\nu(T,x)$ of $T\in \mathcal{T}_{\alpha,\alpha^{\prime}}(\mathbb{F}_{a})$. When $a=0$, ($\mathbb F_a=\mathbb P^1\times \mathbb P^1$), for any current $T\in \mathcal T_{\alpha,\alpha'}(\mathbb P^1\times \mathbb P^1)$, we show that $E^{+}_{(\alpha+\alpha')/3}(T)$ is contained in a curve of total degree $2$, possibly except $1$ point. For any current $T\in \mathcal T_{\alpha,\alpha'}(\mathbb F_a)$, we show that $ E^{+}_{\beta}(T)$ is contained in either in a curve of bidegree $(0,1)$ or in $a+1$  curves of bidegree $(1,0)$ where $\beta\geq (\alpha + (a+1)\alpha^{\prime})/(a+2)$.\\



\end{abstract}

\maketitle

\section{Introduction}

Let $T$ be a positive closed current of bidimension $(p,p)$ on a complex manifold $X$. We consider the upper level sets \[  E_{\alpha}(T)=\{z\in X | \nu(T,z) \geq \alpha \}, \qquad E^{+}_{\alpha}(T)=\{z\in X | \nu(T,z) > \alpha \} \] 
where $\nu(T,z) $ is the Lelong number of $T$ at $z\in X$ and $\alpha \geq 0$ (see \cite{Dem} for the definition of Lelong numbers). It is well known by Siu's result \cite{Siu} that when $\alpha >0$,  $E_{\alpha}(T)$ is an analytic subvariety of $X$ of dimension at most $p$. \\
When $X=\mathbb P^n$ and $T$ is a closed positive current of bidimension $(1,1)$ (or bidegree $(n-1,n-1)$), Coman \cite{Com} proved that  $E^+_{2/3}(T)$ is contained in a complex line and  $|E^+_{1/2}(T)\setminus L |\leq  1$ for some complex line $L$.  In dimension $2$, $|E^+_{2/5}(T)\setminus C |\leq 1$ for some conic $C$ \cite{Com}. In our recent paper  \cite{KY}, we showed that $|E^+_{1/3}(T)\setminus C|\leq 2$ for some cubic curve $C$. In \cite{Com, CT, Hef1, Hef2} further geometric properties of upper level sets $|E^+_{\alpha}(T)|$ are obtained where $T$ is a positive closed current of bidimension $(p,p)$ in $\mathbb P^n$. \\

Geometry of upper level sets were studied for the currents on $X=\mathbb P^m \times \mathbb P^n$ in \cite{CT} and for the currents on multiprojective spaces $X=\mathbb P^{n_1}\times \dots \times \mathbb P^{n_k}$ in \cite{CH}. Proposition 4.1  in \cite{CT} implies that if $T\in \mathcal{T}_{a,b}(\mathbb P^1\times \mathbb P^1)$ (see Section \ref{CH} for the definition of $\mathcal{T}_{a,b}(X)$) , then the upper level set $E^{+}_{(a+b)/2}(T)$ is contained either in a vertical or horizontal fiber of the projections of $\mathbb P^1\times \mathbb P^1$ onto its factors. Since the vertical fibers of $\mathbb P^1\times \mathbb P^1$ are bidegree $(1,0)$ curves and the horizontal fibers are bidegree $(0,1)$ curves on $\mathbb{P}^{1}\times \mathbb{P}^{1}$, this result gives us information about the threshold for the upper level set to be contained in a total degree $1$ curve. Furthermore, this threshold is optimal. 

In this paper, we first  prove a result for thresholds for such upper level sets to be contained in total degree $2$ curves. More precisely, we prove: 

\begin{thm}\label{thm1} Suppose that $T\in \mathcal{T}_{a,b}(\mathbb{P}^{1}\times \mathbb{P}^{1})$. Say $\alpha\geq (a+b)/3$. Then, $E^{+}_{\alpha}(T)$ is contained in a curve of total degree $2$ (i.e. a curve of bidegree $(2,0), (1,1)$ or $(0,2)$), possibly except $1$ point. 
\end{thm} 
 We will prove this theorem in Section \ref{main1}. Example \ref{ex1} in Section \ref{main1} shows that  $\alpha= (a+b)/3$ is sharp in Theorem \ref{thm1}.  

In Section \ref{sec4}, we consider currents on Hirzebruch surfaces $\mathbb F_a=\mathbb{P}(\mathcal{O}\oplus \mathcal{O}(a))$ and prove the  following geometric property  for upper level sets:

\begin{thm} \label{thm2}
Suppose that $\alpha'\geq 0 $, $\alpha+a\alpha'\geq 0 $ and  $T\in \mathcal{T}_{\alpha,\alpha^{\prime}}(\mathbb{F}_{a})$. Say $\beta\geq (\alpha + (a+1)\alpha^{\prime})/(a+2)$. Then, $E^{+}_{\beta}(T)$ is contained in either in a curve of bidegree $(0,1)$ or in $a+1$ fibers of $\varphi$ (i.e. curves of bidegree $(1,0)$). 
\end{thm} 
Fibers of $\varphi$, curves on Hirzebruch surfaces and the currents $T\in \mathcal T_{\alpha,\alpha'}(\mathbb F_a)$  will be defined in the next chapter.  
\begin{rem} 
Notice that if $a=0$, $\mathbb{F}_{a}=\mathbb{P}^{1}\times \mathbb{P}^{1}$, and the result coincides with Coman and Truong's result (Proposition 4.1 in \cite{CT})  about total degree 1 curves on $\mathbb{P}^{1}\times \mathbb{P}^{1}$
\end{rem}
Example \ref{ex2} in Section \ref{sec4} shows that the threshold in the theorem above is sharp when $\alpha \geq \alpha'$. In the case of $\alpha <\alpha'$, it is not clear to us that whether the threshold in Theorem \ref{thm2} is sharp or not.

\section{Preliminaries} 

\subsection{Hirzebruch Surfaces} 

Let us consider rational ruled surfaces, namely ruled surfaces over $\mathbb{P}^{1}$. By Grothendieck's theorem, every rank 2 locally free sheaf on $\mathbb{P}^{1}$ is isomorphic to $\mathcal{O}\oplus \mathcal{O}(a)$ for a unique non-negative integer $a$, and also note that every projective bundle over the projective line arises as the projectivization of a vector bundle (see for instance \cite{GH}, pp. 515-516). The surface $\mathbb{P}(\mathcal{O}\oplus \mathcal{O}(a))$ is called the $a$-th Hirzebruch surface, and it will be denoted by $\mathbb{F}_{a}$. Notice that  $\mathbb{F}_{0}$ is isomorphic to $\mathbb{P}^{1}\times \mathbb{P}^{1}$. Unless $a=1$, these surfaces will not contain any $-1$-curves and will be minimal surfaces. 

Hirzebruch surfaces can also be written in the form of a GIT quotient: 
\[ \mathbb{F}_{a}=(\mathbb{C}^{2}-\{0\})\times (\mathbb{C}^{2}-\{0\})/\sim \] 
where the equivalence relation $\sim$ identifies points in a given orbit of the following $(\mathbb{C}^{*})^2$ action: 
\[ (\lambda,\mu)\cdot (X_{0},X_{1},Y_{0},Y_{1})= (\lambda X_{0}, \lambda X_{1}, \lambda^{-a}\mu Y_{0}, \mu Y_{1}). \] 
The map sending the equivalence class of $(X_{0},X_{1},Y_{0},Y_{1})$ to $[X_{0}:X_{1}]$ is well-defined and this provides us with the ruling map $\varphi: \mathbb{F}_{a}\rightarrow \mathbb{P}^{1}$. Fibers of $\varphi$ are themselves isomorphic to $\mathbb{P}^{1}$, which demonstrates that $\mathbb{F}_{a}$ is a rational ruled surface. For future reference, denote $ \{0\}\times \mathbb{C}^{2}\cup \mathbb{C}^{2}\times \{0\}$ by $Z$ and the  quotient map from $\mathbb{C}^{4}\setminus Z$ to $\mathbb{F}_{a}$ described above by $\pi$.

\subsection{Picard Group and Curves}
Let $F$ be the class of a fiber of $\varphi$ and $H$ be the divisor class corresponding to the invertible sheaf $\mathcal{O}_{\mathbb{F}_{a}}(1)$. Then it is well-known that 
\[ Pic(\mathbb{F}_{a})=\mathbb{Z}F\oplus \mathbb{Z}H. \] 
Furthermore, $F^2=0$, $H^2=a$ and $F\cdot H=1$. If $a>0$, then there is a unique irreducible curve $C$ on $\mathbb{F}_{a}$ with negative self-intersection, with $C^{2}=-a$. The class of $C$ in $Pic(\mathbb{F}_{a})$ (denoted by the same letter) is given by $C=H-aF$. 

One can easily describe curves representing these cohomology classes in terms of the coordinates above: The fibers are clearly movable, and fixing $[X_{0}:X_{1}]$ to be equal to any value $[c_{0}:c_{1}]$ gives us one such fiber. On the other hand for $a>0$, the unique irreducible curve $C$ with $C^{2}=-a$ is not movable, and it is given by $Y_{0}=0$. Let us also remark that a curve of the form $Y_{1}=0$ represents $H$, and it is movable: it can be thought of as coming from the section $(1,0)$ of the vector bundle $\mathcal{O}\oplus \mathcal{O}(a)$. Taking any other section $(1,\sigma)$ of this bundle would give us the moving family (and incidentally shows that $H^2=a$). In contrast, $C$ can be identified with the section $(0,1)$ of the bundle $\mathcal{O}(-a)\oplus \mathcal{O}\cong \mathcal{O}\oplus \mathcal{O}(a)$. The fact that $\mathcal{O}(-a)$ has no non-zero sections indicates that $C$ is not movable. 

Curves on $\mathbb{F}_{a}$ can be described by using the GIT quotient description given above. 
Let $(\alpha,\alpha^{\prime})\in \mathbb{Z}^2$ such that $\alpha^{\prime}\geq 0$ and $\alpha+a \alpha^{\prime}\geq 0$. Then a polynomial $P(X_{0},X_{1},Y_{0},Y_{1})$ is called bihomogeneous of bidegree $(\alpha,\alpha^{\prime})$ if 
\[ P (\lambda X_{0}, \lambda X_{1}, \lambda^{-a}\mu Y_{0}, \mu Y_{1})=\lambda^{\alpha}\mu^{\alpha^{\prime}}P(X_{0},X_{1},Y_{0},Y_{1}). \] 
It is clear that the locus $\{P=0\}$ is a well-defined curve in $\mathbb{F}_{a}$. Conversely, it is a fact that every curve on $\mathbb{F}_{a}$  arises in this way. A curve is said to be of bidegree $(\alpha,\alpha^{\prime})$ if it is the zero-locus of a polynomial of bidegree $(\alpha,\alpha^{\prime})$. A linear polynomial in $X_{0},X_{1}$ is of bidegree $(1,0)$, hence fibers representing $F$ are curves of bidegree $(1,0)$. The polynomial $Y_{0}$ is of bidegree $(-a,1)$, hence $C$ is a curve of bidegree $(-a,1)$. The polynomial $Y_{1}$ is of bidegree $(0,1)$, hence any curve representing $H$ is of bidegree $(0,1)$.

\subsection{Currents on Hirzebruch surfaces} \label{CH}

Following Guedj \cite{G}, we now describe currents on $\mathbb{F}_{a}$. Let $\mathcal{T}(\mathbb{F}_{a})$ denote the cone of positive, closed currents of bidegree $(1,1)$ on $\mathbb{F}_{a}$. If $T\in \mathcal{T}(\mathbb{F}_{a})$ then $\pi^{*}(T)$ is a positive, closed current of bidegree $(1,1)$ on $\mathbb{C}^{4}\setminus Z$. However, $Z$ is of codimension $2$ in $\mathbb{C}^{4}$, hence $\pi^{*}(T)$ extends to a $(1,1)$-current on $\mathbb{C}^{4}$. By the $\partial \overline{\partial}$-Poincar\'e lemma for currents, there exists a plurisubharmonic function $\psi$ on $\mathbb{C}^{4}$ such that $\pi^{*}T=dd^{c}\psi$. Theorem 3.1 in \cite{G} strengthens this as follows: Let $(\alpha,\alpha^{\prime})\in \mathbb{R}^2$. Define $\mathcal{P}_{\alpha,\alpha^{\prime}}$ to be the set of PSH-functions $\psi$ on $\mathbb{C}^{4}$ having their supremum on the unit ball in $\mathbb{C}^4$ equal to $0$ and such that 
\[ \psi(\lambda X_{0}, \lambda X_{1},  \lambda^{-a}\mu Y_{0}, \mu Y_{1})= \alpha \log|\lambda|+ \alpha^{\prime} \log |\mu|+\psi(X_{0},X_{1},Y_{0},Y_{1}). \]
The cone of positive, closed currents of bidegree $(1,1)$ on $\mathbb F_a$ decomposes as $\mathcal T(\mathbb F_a)=\cup \mathcal T_{\alpha,\alpha'}(\mathbb F_a)$  where $\mathcal{T}_{\alpha,\alpha^{\prime}}(\mathbb{F}_{a})$ denotes the set of positive, closed $(1,1)$-currents on $\mathbb{F}_{a}$ whose cohomology class is $\alpha F+\alpha^{\prime} H$.  We should note that the class of $\alpha F+\alpha^{\prime} H$ is pseudoeffective,  i.e.  $\mathcal{T}_{\alpha,\alpha^{\prime}}(\mathbb{F}_{a})\neq \emptyset $ if and only if $\alpha'\geq 0, \alpha+a\alpha'\geq 0.$    Guedj shows that for $T\in \mathcal{T}_{\alpha,\alpha^{\prime}}(\mathbb{F}_{a})$ there exists unique $\psi\in \mathcal{P}_{\alpha,\alpha^{\prime}}$ such that $\pi^{*}(T)=dd^{c} \psi$. 

The currents on $\mathbb F_0=\mathbb P_1^1 \times \mathbb P_2^1$ can be written in terms of the Fubini-Study form on $\mathbb P^1$. Let $$\pi_1: \mathbb P_1^1 \times \mathbb P_2^1 \to \mathbb P_1^1, \,  \pi_2: \mathbb P_1^1 \times \mathbb P_2^1 \to \mathbb P_2^1 $$  be the canonical projections. Define $$\omega_1=\pi_1^*\omega, \, \omega_2=\pi_2^*\omega $$ where $\omega$ is the Fubini-Study form on $\mathbb P^1$.  Then $\mathcal{T}_{a,b}(\mathbb P^1 \times \mathbb P^1)$ is the set of positive, closed $(1,1)$-currents on $\mathbb P^1 \times \mathbb P^1$ whose cohomology class is $a\omega_1 + b \omega_2$.

\section{Thresholds for Lelong Numbers on $\mathbb{P}^{1}\times \mathbb{P}^{1}$ for Degree 2 Curves} \label{main1}

Let us consider $X=\mathbb{P}^{1}\times \mathbb{P}^{1}$ in this section. Suppose that $T\in \mathcal{T}_{a,b}$. Let $\nu(T,z)$ denote the Lelong number of the current $T$ at a point $z\in \mathbb{P}^{1}\times \mathbb{P}^{1}$. For $\alpha\geq 0$, define the upper level sets 
\[  E_{\alpha}(T)=\{z\in \mathbb{P}^{1}\times \mathbb{P}^{1} | \nu(T,z) \geq \alpha \}, \qquad E^{+}_{\alpha}(T)=\{z\in \mathbb{P}^{1}\times \mathbb{P}^{1} | \nu(T,z) > \alpha \}.   \] 
By Siu's theorem \cite{Siu}, we know that $E_{\alpha}(T)$ is an analytic subvariety of $\mathbb{P}^{1}\times \mathbb{P}^{1}$ of dimension at most $1$. Let $\pi_{1}$ and $\pi_{2}$ denote the projection maps from $\mathbb{P}^{1}\times \mathbb{P}^{1}$ to its factors. Let us first show that the threshold for the upper level sets in Theorem \ref{thm1} is sharp.

\begin{ex} \label{ex1}
Let $V_{1}, V_{2}, V_{3}$ be three distinct fibers of $\pi_{1}$ and $H_{1}, H_{2}, H_{3}$ be three distinct fibers of $\pi_{2}$. Denote the currents of integration along these curves by the same letters. Then the $(1,1)$-current 
\[ T=\frac{a}{3}(V_{1}+V_{2}+V_{3})+\frac{b}{3}(H_{1}+H_{2}+H_{3}) \] 
lies in $\mathcal{T}_{a,b}$. Let $S$ denote the set of nine points of intersections of $V_{i}$'s with $H_{j}$'s. Then for any $p\in S$, we have $\nu(T,p)=(a+b)/3$. Therefore, if $\alpha< (a+b)/3$, then $S\subset E^{+}_{\alpha}(T)$. We claim that $E^{+}_{\alpha}(T)$ cannot be contained in a curve of total degree 2 possibly except $1$-point: First, let us consider the bidegree $(1,1)$ case. Such a curve would be either irreducible or reducible. An irreducible $(1,1)$-curve intersects any horizontal or vertical fiber at $1$ point. Therefore, it can contain at most $3$ points from $S$. A reducible $(1,1)$-curve is the union of a vertical fiber and a horizontal fiber, hence it can contain at most $5$ points from $S$. Then let us consider the cases of bidegree $(2,0)$ or $(0,2)$. Such curves are necessarily reducible, and can contain at most $6$ points of $S$. In any case, there are at least $3$ points of $S$ not contained in the curve, which proves our claim. This example shows that the threshold $(a+b)/3$  in Theorem \ref{thm1} is optimal. 
\end{ex} 

\begin{proof}(of Theorem \ref{thm1}) 
Assume, to the contrary, that there doesn't exist any curve of total degree $2$ containing $E^{+}_{\alpha}(T)$  possibly except $1$ point. Note that a $(1,1)$-curve on $\mathbb{P}^{1}\times \mathbb{P}^{1}$ is the zero locus of a polynomial of the form 
$a_{00}X_{0}Y_{0}+a_{01}X_{0}Y_{1}+a_{10}X_{1}Y_{0}+a_{11}X_{1}Y_{1}$, so there exists a $(1,1)$-curve passing through any $3$ points in $\mathbb{P}^{1}\times \mathbb{P}^{1}$ (this curve is unique unless all $3$ points lie on a vertical fiber or on a horizontal fiber). Therefore, the assumption implies that $E^{+}_{\alpha}(T)$ must contain a subset $A=\{P_{1},P_{2},P_{3},P_{4},P_{5}\}$ of $5$ points. Let $m_{i,j}(A)$ denote the maximum number of points of $A$ contained in a curve of bidegree $(i,j)$. 

\textbf{Case 1:} Suppose that $m_{1,1}(A)=3$. In this case, neither of $m_{1,0}(A)$ and $m_{0,1}(A)$ can be more than $2$. 
Let us first consider the subcase where $m_{1,0}(A)=m_{0,1}(A)=1$.  Let $\mathcal{X}$ be the space of bihomogeneous polynomials of bidegree $(2,2)$ passing through all points of $A$ and with multiplicity at least $2$ at $P_{1}$. Passing through each $P_{i}$ imposes 1 linear condition and the vanishing of two partial derivatives at $P_{1}$ imposes 2 more linear conditions, hence $dim(\mathcal{X})\geq 9-(5+2)=2$. Take a $(1,1)$-curve $C_{1}$ through $P_{1}, P_{2}, P_{3}$ and a $(1,1)$-curve through $C_{2}$ through $P_{1}, P_{4}, P_{5}$. (Under the current assumptions, these curves are uniquely determined and each contain precisely $3$ points from $A$.) Then, clearly $Q_{1}=C_{1}C_{2}\in \mathcal{X}$, where we allow the ambiguity in the notation so as to denote the equation of curves with the same letters. Since $dim(\mathcal{X})\geq 2$, there exists $Q_{2}\in \mathcal{X}$ such that $Q_{1}, Q_{2}$ are linearly independent. We claim that $Q_{1}$ and $Q_{2}$ do not have a non-trivial common factor. Indeed, first observe that $C_{1}$ and $C_{2}$ have to be irreducible since 
 $m_{1,0}(A)=m_{0,1}(A)=1$. Since irreducible $(1,1)$-curves in $\mathbb{P}^1\times \mathbb{P}^1$ are smooth, $C_{1}$ and $C_{2}$ have to be smooth. If $C_{1}$ divides $Q_{2}$, then $Q_{2}=C_{1}C$ for some degree $(1,1)$-curve. Since $C_{1}$ doesn't pass through $P_{4}, P_{5}$ and has multiplicity $1$ at $P_{1}$ (given that it is smooth), the curve $C$ must pass through $P_{1}, P_{4}$ and $P_{5}$. But now, by Bezout's theorem, the intersection number of $C$ and $C_{2}$ being greater than $2$ implies that they share a component, hence $C=C_{2}$. This contradicts the linear independence of $Q_{1}$ and $Q_{2}$. A similar argument shows that $C_{2}$ cannot divide $Q_{2}$, proving the claim. Now, define 
\[  u=\frac{1}{2} \log(|Q_{1}|^2+|Q_{2}|^2).  \]
Then by Theorem 3.1 in \cite{G}, $dd^{c} u$ determines a  current  $S\in \mathcal T_{2,2}(\mathbb P^1\times \mathbb P^1)$. Since $S$  is smooth except at finitely many logarithmic poles,  
by a result of Demailly \cite{Dem},  the intersection $T\wedge S$ is a well-defined positive measure. We observe that 
\[ 2(a+b)=\int_X T\wedge S \geq   \sum_{i=1}^{5} T\wedge S (P_i)\geq  \sum_{i=1}^{5} \nu(T,P_{i})\nu(S,P_{i})> 6\alpha. \] 
The second inequality follows from Demailly's comparison theorem for Lelong numbers  (Corollary 5.10 in \cite{Dem}).  The last inequality follows from $\nu(T,P_{i})>\alpha$ for all $i$, $\nu(S,P_{1})\geq 2$ and $\nu(S,P_{j})\geq 1$ for $j=2,3,4,5$. We deduce that $\alpha<(a+b)/3$, which is a contradiction. 

The remaining subcase is when at least one of $m_{1,0}(A)$ or $m_{0,1}(A)$ is equal to $2$. Assume that $m_{0,1}(A)=2$ without loss of generality. We will refer to $(0,1)$-curves as horizontal lines and $(1,0)$-curves as vertical lines. 
Assume, without loss of generality again, that $P_{1}$ and $P_{2}$ lie on a horizontal line. First suppose that there exists a point of $A$, say $P_{5}$, which does not lie on a horizontal or vertical line through any other point of $A$ (see Figure \ref{fig1}). We may  assume that $P_{3}$ does not lie on the vertical line through $P_{2}$ and $P_{4}$ does not lie on the vertical line through $P_{1}$, if necessary by reordering $P_{3}$ and $P_{4}$. Now, the $(1,1)$-curve $C_{1}$ through $P_{1}, P_{4}, P_{5}$ and the $(1,1)$-curve $C_{2}$ through $P_{2}, P_{3},P_{5}$ both have to be irreducible. Therefore we can repeat the argument for the subcase $m_{1,0}(A)=m_{0,1}(A)=1$ to finish the proof (with $P_{5}$ instead of $P_{1}$, etc.). The remaining case is when every point of $A$ is on a vertical or a horizontal line through another point of $A$. But in this case there must exist a vertical and a horizontal line which contain in total $4$ points of $A$ (see Figure \ref{fig2}). This contradicts the assumption that $m_{1,1}(A)=3$. 
\begin{figure}[h]
\begin{center} 
\scalebox{0.35}
{\includegraphics{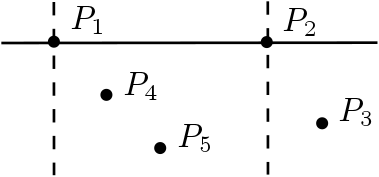}}
\caption{A configuration of $A$ when $m_{1,1}(A)=3, m_{0,1}(A)=2$.}\label{fig1} 
\end{center}
\end{figure}

\begin{figure}[h]
\begin{center} 
\scalebox{0.20}
{\includegraphics{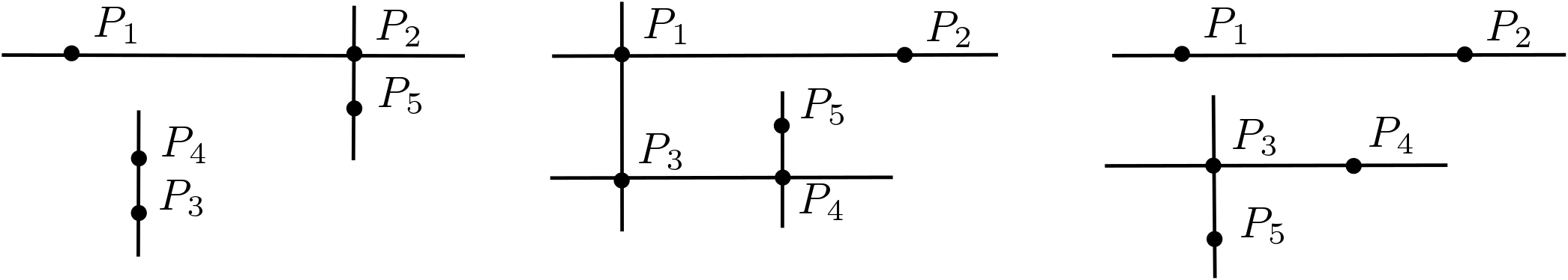}}
\caption{Other configurations of $A$ when $m_{1,1}(A)=3, m_{0,1}(A)=2$. }\label{fig2} 
\end{center}
\end{figure}

\textbf{Case 2:} Suppose that $m_{1,1}(A)=4$ and $m_{1,0}(A)=m_{0,1}(A)=1$. Say $C_{1}$ is a curve of bidegree $(1,1)$ passing through $P_{1},P_{2},P_{3}, P_{4}$ without loss of generality. Clearly, $C_{1}$ must be irreducible and $P_{5}\notin C_{1}$. Since $E^{+}_{\alpha}(T)$ is not contained in a bidegree $(1,1)$-curve except one point, by the standing assumption, there exists a sixth point $P_{6}\in E^{+}_{\alpha}(T)$ such that $P_{6}\notin C_{1}$. Suppose first that $P_{6}$ is not on a common vertical or horizontal line with $P_{5}$. Since $P_{6}$ can be on a common vertical (resp. horizontal) line with at most 1 point from $\{P_{1},P_{2},P_{3},P_{4}\}$, we may assume without loss of generality that $P_{6}$ is not on a common vertical or horizontal line with $P_{1}$(see Figure \ref{fig3}(a)). Then, the $(1,1)$-curve $C_{2}$ passing through $P_{1}, P_{5}$ and $P_{6}$ must be irreducible.  The intersection number of $C_{1}$ and $C_{2}$ is $2$, hence $C_{2}$ can  contain at most one more point from $A$, say $P_{4}$. Now, we can conclude by repeating the argument in the first part of Case 1 for the irreducible curves $C_{1}, C_{2}$ and the points $P_{1}, P_{2}, P_{3}, P_{5}, P_{6}$. 

\begin{figure}[h]
\begin{center} 
\scalebox{0.40}
{\includegraphics{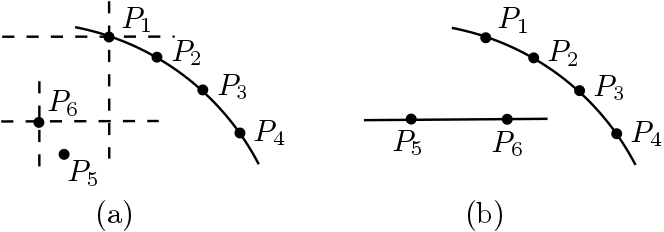}}
\caption{Configurations of $A$ when $m_{1,1}(A)=4, m_{0,1}(A)=m_{1,0}(A)=1$.}\label{fig3} 
\end{center}
\end{figure}

The remaining subcase is where $P_{5}$ and $P_{6}$ lie on a common vertical or horizontal line (see Figure \ref{fig3}(b)). Since these two cases are similar, let us assume without loss of generality that they lie on a common horizontal line. Let $B=\{P_{1}, P_{2}, P_{3}, P_{5}, P_{6}\}$. We claim that $m_{1,1}(B)=3$: Indeed, suppose that there exists a $(1,1)$-curve $C$ passing through 4 (or more) points of $B$. If $C$ contains all of $P_{1}, P_{2}, P_{3}$ then its intersection number with $C_{1}$ exceeds $2$. But $C_{1}$ is irreducible, so $C=C_{1}$. This is a contradiction since neither of $P_{5}, P_{6}$ lie on $C_{1}$. The remaining case is when $C$ contains both $P_{5}$ and $P_{6}$, and two of $P_{1}, P_{2}, P_{3}$. But then $C$ must be reducible as the union of the horizontal line through $P_{5},P_{6}$ and a vertical line through two of $P_{1}, P_{2}, P_{3}$. This contradicts $m_{1,0}(A)=1$, proving the claim. Replacing the set $A$ by the set $B$, we go back to Case 1. 

\textbf{Case 3:}  Suppose that $m_{1,1}(A)=4$, $m_{1,0}(A)=2$, and $m_{0,1}(A)\leq 2$. As the first subcase, assume that there exists an irreducible $(1,1)$-curve $C_{1}$ passing through $4$ of the points in $A$, say through $P_{1}, P_{2}, P_{3}, P_{4}$. Since $m_{1,0}(A)=2$, there exists a vertical line containing 2 points from $A$. One of these points must be $P_{5}$ since $C_{1}$ is irreducible. Without loss of generality, assume that the other point on this vertical line is $P_{1}$. Then, $P_{5}$ can lie on the same horizontal line with at most $1$ point in $A$, say $P_{4}$. Now, the $(1,1)$-curve $C_{2}$ passing through $P_{2}, P_{3}$ and $P_{5}$ is irreducible. One can then repeat the argument in the first subcase of Case 1 with these two curves to conclude ($P_{3}$ here will play the role of $P_{1}$ in that subcase). 

%The fact that $C_{1}$ passes also through $P_{2}$ in addition to $P_{1}, P_{3}, P_{4}$ makes the final inequality stronger, giving us $7\alpha$ instead of $6\alpha$). 

The remaining case is when there is no irreducible $(1,1)$-curve through any $4$ points of $A$. Then a $(1,1)$-curve through $4$ points of $A$, which exists by assumption, must be a union of a horizontal line and a vertical line. Since $m_{1,0}(A)$ and $m_{0,1}(A)$ are both at most $2$, both the horizontal and the vertical line must contain two points each. Without loss of generality, suppose that $P_{1}, P_{2}$ lie on a horizontal line and $P_{3}, P_{4}$ lie on a vertical line.  Also, we may assume without loss of generality that $P_{5}$ and $P_{3}$ are not on the same horizontal line, and $P_{5}$ and $P_{1}$ are not on the same vertical line. Then, the $(1,1)$-curve $C_{1}$ passing through $P_{1}, P_{3}, P_{5}$ is irreducible. If also $P_{5}$ does not lie on the vertical line through $P_{2}$ and does not lie on the horizontal line through $P_{4}$ (see Figure \ref{fig4}) , then the $(1,1)$-curve $C_{2}$ passing through $P_{2}, P_{4}, P_{5}$ is also irreducible, and we can use $C_{1}, C_{2}$ as in the first subcase of Case 1. \\
\begin{figure}[h]
\begin{center} 
\scalebox{0.12}
{\includegraphics{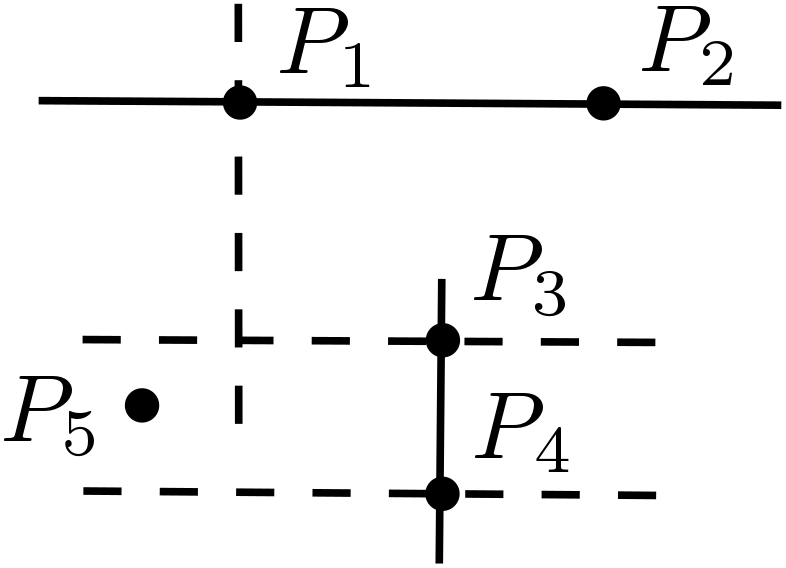}}
\caption{A configuration of $A$ when $m_{1,1}(A)=4, m_{1,0}(A)=m_{0,1}(A)=2$.}\label{fig4} 
\end{center}
\end{figure}

Let us now exhaust the remaining possibilities. We may assume that $P_{4}$ and $P_{5}$ lie on the same horizontal line without loss of generality.  Since $E^{+}_{\alpha}(A)$ is not contained in a $(1,1)$-curve except one point by the standing assumption, there must exist another point $P_{6}\in E^{+}_{\alpha}(A)$ which does not lie on the first $(1,1)$-curve, namely does not lie on either of the lines through $P_{1}, P_{2}$ or $P_{3}, P_{4}$. We finally analyze the possibilities for the position of $P_{6}$: If $P_{6}$ is not on the horizontal line through $P_{4}$ and $P_{5}$ and also $P_{2}$ and $P_{6}$ are not on a vertical fiber(see Figure \ref{fig5}(a)), then the curve $C_{2}$ through $P_{2}, P_{4}, P_{6}$ is irreducible. Using the curves $C_{1}, C_{2}$ as in the first subcase of Case 1 will finish the argument. Likewise, if  $P_{6}$ is not on the horizontal line through $P_{4}$ and $P_{5}$ and also $P_{1}$ and $P_{6}$ are not on a vertical fiber, then the curve $C_{2}$ through $P_{1}, P_{4}, P_{6}$ is irreducible and again $C_{1}, C_{2}$ can be used in the same way. The only remaining case is when $P_{6}$ lies on the horizontal line through $P_{4}$ and $P_{5}$. If $P_{2}$ and $P_{6}$ are not on the same vertical line (See Figure \ref{fig5}(b)), then the curve $C_{2}$ through $P_{2}, P_{3}, P_{6}$ is irreducible and we can finish as before. If $P_{2}$ and $P_{6}$ are on the same vertical line (see Figure \ref{fig5}(c)), then set $B=\{P_{1}, P_{2}, P_{3}, P_{5}, P_{6}\}$. We claim that $m_{1,1}(B)=3$. Indeed, a $(1,1)$-curve passing through $4$ points of $B$ must contain at least one of $P_{2}$ or $P_{6}$. But both of these points are on junctions of vertical and horizontal lines each containing $2$ points of $B$, so this implies that such a $(1,1)$-curve would have to be reducible. But then we immediately check that no reducible $(1,1)$-curve contains $4$ points of $B$, which proves the claim. Replacing $A$ with $B$ takes us back to Case 1.  
 \begin{figure}[h]
\begin{center} 
\scalebox{0.45}
{\includegraphics{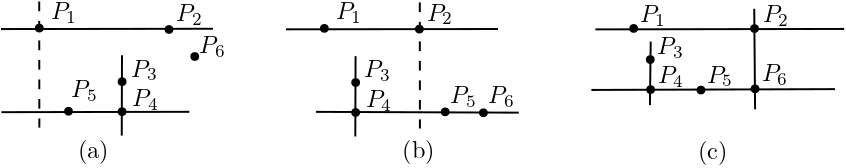}}
\caption{Other configurations of $A\cup \{P_6\}$ when $m_{1,1}(A)=4, m_{1,0}(A)=m_{0,1}(A)=2$.}\label{fig5} 
\end{center}
\end{figure}
 
\textbf{Case 4:} Suppose that $m_{1,1}(A)=4$, $m_{1,0}(A)=2$ and $m_{0,1}(A)=3$. Since a horizontal line and a vertical line together cannot contain $5$ distinct points in view of $m_{1,1}(A)=4$, we may assume without loss of generality that $P_{1}, P_{2}, P_{3}$ lie on a horizontal line, and $P_{1}, P_{4}$ lie on a vertical line. Then $P_{5}$ does not lie on these two lines and we may also assume without loss of generality that $P_{5}$ does not lie on the same vertical line with $P_{2}$. Also, there must exist $P_{6}\in  E^{+}_{\alpha}(A)$ which does not lie on the first two lines (see Figure \ref{fig6}). Let $B=\{P_{1}, P_{2}, P_{4}, P_{5}, P_{6}\}$. If $m_{1,1}(B)=3$, then we can replace $A$ with $B$ and go back to Case 1. Otherwise, we must have $m_{1,1}(B)=4$. It is clear that no $3$ points of $B$ are on the same vertical line, so $m_{1,0}(B)=2$. If $m_{0,1}(B)=2$ as well, then replace $A$ with $B$ and return to Case 3. 

\begin{figure}[h]
\begin{center} 
\scalebox{0.30}
{\includegraphics{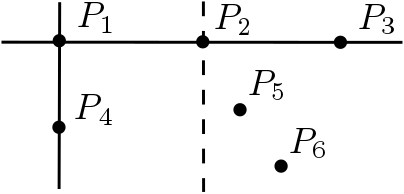}}
\caption{Configuration of $A\cup \{P_6\}$ when $m_{1,1}(A)=4, m_{1,0}(A)=2,m_{0,1}(A)=3$.}\label{fig6} 
\end{center}
\end{figure}

The only remaining possibillity is that $m_{0,1}(B)=3$. This can happen only if $P_{4}, P_{5}, P_{6}$ are on the same horizontal line. Since $E^{+}_{\alpha}(A)$ is not contained in a degree $(0,2)$-curve possibly except one point, there must exist $P_{7},P_{8}\in E^{+}_{\alpha}(A)$ which do not lie on the two horizontal lines through the previous 6 points (see Figure \ref{fig7}). If $P_{7}$ does not lie on the vertical line through $P_{1}$ and $P_{4}$, then set $D=\{P_{1}, P_{2}, P_{4}, P_{5}, P_{7}\}$. We then have $m_{1,1}(D)=3$ since any irreducible $(1,1)$-curve can pass through at most one of $P_{1}, P_{2}$ and at most one of $P_{4}, P_{5}$, furthermore any reducible $(1,1)$-curve can pass through at most $3$ of these points by the given assumptions. Hence, we can replace $A$ by $D$ and go back to Case 1. If $P_{7}$ lies on the vertical line through $P_{1}$ and $P_{4}$, then let $D=\{P_{2}, P_{3}, P_{4}, P_{5}, P_{7}\}$. By the choices made, we have $m_{1,1}(D)=4$ and $m_{1,0}(D)=m_{0,1}(D)=2$, so replacing $A$ by $D$ we go back to Case $3$. This finishes the proof in Case 4. 

\begin{figure}[h]
\begin{center} 
\scalebox{0.18}
{\includegraphics{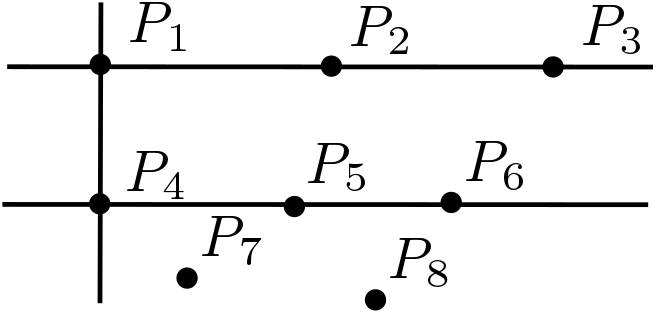}}
\caption{Special configuration of $A\cup \{P_6\}$ when $m_{1,1}(A)=4, m_{1,0}(A)=2,m_{0,1}(A)=3$.}\label{fig7} 
\end{center}
\end{figure}

\textbf{Case 5:} Suppose that $m_{1,1}(A)=5$. First assume that there exists an irreducible $(1,1)$-curve $C$ passing through all points of $A$. Then by the assumption, $E^{+}_{\alpha}(A)$  is not contained in $C$ possibly except $1$-point. Hence, there must exist $P_{6}, P_{7}\in E^{+}_{\alpha}(A)$ which are not on $C$. Let $B=\{P_{1}, P_{2}, P_{3}, P_{6}, P_{7}\}$. Then $m_{1,1}(B)\leq 4$ since any $(1,1)$-curve other than $C$ intersects $C$ in at most $2$ points, hence misses at least 1 point of $B$. Replacing $A$ by $B$, we go back to one of the cases 1, 2, 3 or 4. 

The remaining subcase is when a reducible $(1,1)$-curve passes through all points of $A$. Without loss of generality, we may assume that a horizontal line contains at least as many points of $A$ as any vertical line.

\begin{figure}[h]
\begin{center} 
\scalebox{0.30}
{\includegraphics{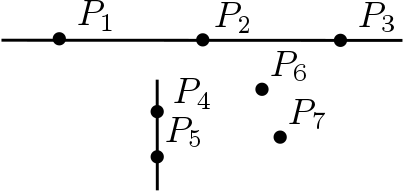}}
\caption{The configuration of $A$ when $m_{1,1}(A)=5, m_{0,1}(A)=3$.}\label{fig8} 
\end{center}
\end{figure}

Suppose first that $P_{1}, P_{2}, P_{3}$ lie on a horizontal line which doesn't contain $P_{4}, P_{5}$, and  so $P_{4}, P_{5}$ lie on a vertical line. Then, there must exist $P_{6}, P_{7}\in E^{+}_{\alpha}(A)$ which are not on these two lines (see Figure \ref{fig8}). Set $B=\{P_{2}, P_{3}, P_{4}, P_{6}, P_{7}\}$. Then, no $(1,1)$-curve can contain all points of $B$: An irreducible curve would miss at least one of $P_{2}, P_{3}$ and a reducible curve can't contain them all by the choices made. Hence $m_{1,1}(B)\leq 4$ and we can replace $A$ by $B$ in order to go back to one of cases 1, 2, 3 or 4. 

\begin{figure}[h]
\begin{center} 
\scalebox{0.26}
{\includegraphics{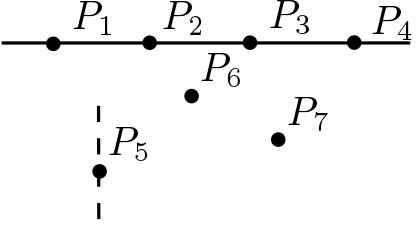}}
\caption{The configuration of $A$ when $m_{1,1}(A)=5, m_{0,1}(A)=4$.}\label{fig9} 
\end{center}
\end{figure}

Second, assume that $P_{1}, P_{2}, P_{3}, P_{4}$ are on a horizontal line which does not contain $P_{5}$. Again,  there must exist $P_{6}, P_{7}\in E^{+}_{\alpha}(A)$ which are not on this horizontal line and the vertical line through $P_{5}$ (see Figure \ref{fig9}). Set $B=\{P_{3}, P_{4}, P_{5}, P_{6}, P_{7}\}$. Again, a $(1,1)$-curve cannot contain all points of $B$: An irreducible $(1,1)$-curve would miss at least one of $P_{3}$ and $P_{4}$, and a reducible curve cannot contain them all by the choices made. Hence $m_{1,1}(B)\leq 4$ and we can replace $A$ by $B$ in order to go back to one of cases 1, 2, 3 or 4 again. 

Finally, assume that $P_{1}, P_{2}, P_{3}, P_{4}, P_{5}$ lie on a horizontal line. Then  there must exist $P_{6}, P_{7}\in E^{+}_{\alpha}(A)$ which are not on this horizontal line and which don't lie on a common vertical line. Set $B=\{P_{3}, P_{4}, P_{5}, P_{6}, P_{7}\}$. Similar to the subcases above, no $(1,1)$-curve can contain all points of $B$: An irreducible curve would miss at least two of $P_{3}, P_{4}, P_{5}$ and a reducible curve cannot contain them all either. Hence $m_{1,1}(B)\leq 4$ and we can replace $A$ by $B$ in order to go back to one of cases 1, 2, 3 or 4. This finishes the proof of this case and also the proof of the theorem. 
\end{proof} 

\section{Thresholds for Lelong Numbers on $\mathbb{F}_{a}$} \label{sec4}

In this section we consider currents on Hirzebruch surfaces $\mathbb{F}_{a}$. As defined in the introduction, let  $\mathcal{T}_{\alpha,\alpha^{\prime}}(\mathbb{F}_{a})$ denote the set of positive, closed $(1,1)$-currents on $\mathbb{F}_{a}$ whose cohomology class is $\alpha F+\alpha^{\prime} H$. Let $T\in \mathcal{T}_{\alpha,\alpha^{\prime}}(\mathbb{F}_{a})$. For $\beta\geq 0$ define 
\[  E_{\beta}(T)=\{z\in \mathbb{F}_{a} | \nu(T,z) \geq \beta \}, \qquad E^{+}_{\beta}(T)=\{z\in \mathbb{F}_{a} | \nu(T,z) > \beta \}.   \] 
By Siu's theorem, for any $\beta>0$,  $E_{\beta}(T)$ is an analytic subvariety of $\mathbb{F}_{a}$ of dimension at most $1$. Let us first show that the threshold for the upper level sets in Theorem \ref{thm2} is sharp when $\alpha\geq \alpha'$. 

\begin{ex} \label{ex2}
Let $C$ be the (unique) curve of bidegree $(-a,1)$ on $\mathbb{F}_{a}$ and let $F_{1}, F_{2},\ldots, F_{a+2}$ be $a+2$ distinct fibers of $\varphi$. Denote the currents of integration along these curves by the same letters. Then the $(1,1)$-current 
\[ T= \left( \frac{\alpha+(a+1)\alpha^{\prime}}{a+2}\right) (F_{1}+F_{2}+\ldots + F_{a+1})+ \left( \frac{\alpha-\alpha^{\prime}}{a+2}\right) F_{a+2} +\alpha^{\prime} C \] 
lies in $\mathcal{T}_{\alpha, \alpha^{\prime}}(\mathbb{F}_{a})$ when $\alpha \geq \alpha'$. Let $P$ be the point of intersection of $F_{a+2}$ with $C$. Let 
\[ S= F_{1}\cup F_{2}\cup \ldots \cup F_{a+1} \cup \{P\}. \] 
Then for every $q\in S$ we have $\nu(T,q)= (\alpha + (a+1)\alpha^{\prime})/(a+2)$. Suppose we now choose $\beta$ such that $\beta<   (\alpha + (a+1)\alpha^{\prime})/(a+2)$. Then we will have $S\subset E^{+}_{\beta}(T)$. The set $S$ is clearly not contained in $a+1$ fibers of $\varphi$. We claim that it is not contained in a curve of bidegree $(0,1)$ either: An irreducible $(0,1)$-curve can only contain finitely many points from each fiber of $\varphi$, hence cannot contain $S$. A reducible $(0,1)$-curve has to be the union of $C$ and $a$ fibers of $\varphi$, hence it cannot contain the $a+1$ distinct  fibers. This shows that the bound for $\beta$ in Theorem \ref{thm2} is sharp when $\alpha\geq \alpha'$. 
\end{ex} 

\begin{proof}(of Theorem \ref{thm2}) 
Let us first note that a $(0,1)$-curve on $\mathbb{F}_{a}$ is the zero locus of a polynomial which is a linear combination of the monomials $Y_{1}, X_{0}^{a}Y_{0}, X_{0}^{a-1}X_{1}Y_{0}, \ldots, X_{1}^{a}Y_{0}$. In particular, there exists a $(0,1)$-curve passing through any $a+1$ points in $\mathbb{F}_{a}$. Suppose, to the contrary to the statement of the theorem, that $E_{\beta}^{+}(T)$ is not contained in a curve of bidegree $(0,1)$ or $a+1$ fibers of $\varphi$. Then, in particular, there must exist a subset $A=\{Q_{1}, Q_{2}, \ldots, Q_{a+2}\}$ of $E_{\beta}^{+}(T)$  such that no two of $Q_{i}$'s lie on the same fiber of $\varphi$.  

Let $\mathcal{Y}$ be the vector space of bihomogeneous polynomials of bidegree $(1,1)$, namely the set of $P(X_{0},X_{1}, Y_{0}, Y_{1})$ such that
\[ P(\lambda X_{0}, \lambda X_{1}, \lambda^{-a}\mu Y_{0}, \mu Y_{1})= \lambda \mu P(X_{0},X_{1},Y_{0}, Y_{1}). \] 
Since a basis for $\mathcal{Y}$ is given by the monomials $X_{0}Y_{1}, X_{1}Y_{1}, X_{0}^{a+1}Y_{0}, X_{0}^{a}X_{1}Y_{0}, \ldots, X_{1}^{a+1}Y_{0}$, we observe that $\dim(\mathcal{Y})=a+4$. Let $\mathcal{X}$ be the subspace of $\mathcal{Y}$ vanishing at all points of $A$. Since each point of $A$ imposes at most one independent linear condition, we get that $\dim(\mathcal{X})\geq 2$. Let $C_{1}$ be a curve of bidegree $(0,1)$ through $Q_{1},\ldots, Q_{a+1}$ and $C_{2}$ be a curve of bidegree $(1,0)$ through $Q_{a+2}$. Then $C_{2}$ doesn't pass through any $Q_{j}$ with $j\leq a+1$ by our initial assumption. We have $P_{1}=C_{1}C_{2}\in \mathcal{X}$ and since $\dim(\mathcal{X})\geq 2$, there exists $P_{2}\in \mathcal{X}$ such that $P_{1}$ and $P_{2}$ are linearly independent. 

We now want to show that $P_{1}$ and $P_{2}$ above do not have any non-trivial common factor, if necessary after modifying the choice of $A$. There are several subcases: 

\textbf{Case 1:} First, let us assume that $C_{1}$ is irreducible and $m_{0,1}(A)=a+1$, that is  the maximum number of points of $A$ contained in a curve of bidegree $(0,1)$ is $a+1$. If $C_{2}$ divides $P_{2}$, then we can write $P_{2}=C  C_{2}$ where $C$ is of bidegree $(0,1)$. Then, $C$ would have to contain $Q_{1},\ldots,Q_{a+1}$ since $C_{2}$ doesn't contain any of these points. But now, $C$ and $C_{1}$ intersect at $a+1$ points, which is 1 more than the upper bound allowed by Bezout's theorem. Hence, $C=C_{1}$. This contradicts the linear independence of $P_{1}$ and $P_{2}$. If $C_{1}$ divides $P_{2}$, then we can write $P_{2}=C C_{1}$ where $C$ is of bidegree $(1,0)$. Since $m_{0,1}(A)=a+1$, the point $Q_{a+2}$ cannot be on $C_{1}$. Thus $Q_{a+2}\in C$. This shows that $C=C_{2}$, again contradicting the linear independence of $P_{1}$ and $P_{2}$. 

\textbf{Case 2:} Assume that $C_{1}$ is reducible and $m_{0,1}(A)=a+1$. Then we can write $C_{1}=C^{\prime} F_{1} F_{2}\ldots F_{a}$ where $C^{\prime}$ is the (unique, irreducible) curve of bidegree $(-a,1)$ and $F_{i}$'s are fibers of $\varphi$. The curve $C^{\prime}$ contains exactly $1$ point from $A\setminus \{Q_{a+2}\}$, say $Q_{a+1}$, since otherwise either two points from $A$ would have to be contained in the same fiber of $\varphi$, or $C^{\prime}$ together with $a$ other fibers would contain all $a+2$ points of $A$, violating the assumption $m_{0,1}(A)=a+1$ (see Figure \ref{fig10}). Then, each $F_{i}$ has to contain precisely $1$ point of $A\setminus \{Q_{a+2}\}$, say $Q_{i}$.  Now, if $C^{\prime}$ divides $P_{2}$, then $P_{2}=C^{\prime} C$ where $C$ is a curve of bidegree $(a+1,0)$ passing through $Q_{1}, Q_{2}, \ldots, Q_{a}$ and $Q_{a+2}$. But then we necessarily have $C=F_{1}F_{2}\ldots F_{a} C_{2}$, contradicting the linear independence of $P_{1}$ and $P_{2}$. If some $F_{i}$, say without loss of generality $F_{1}$, divides $P_{2}$, then write $P_{2}=F_{1}C$ where $C$ is a bidegree $(0,1)$ curve passing through $Q_{2}, Q_{3},\ldots Q_{a+2}$. The intersection number of a curve of bidegree $(0,1)$ and the curve of bidegree $(-a,1)$ on $\mathbb{F}_{a}$ is $0$. However, $C$ and $C^{\prime}$ share the point $Q_{a+1}$. This implies that $C$ is divisible by $C^{\prime}$. Writing $C=C^{\prime}C^{\prime \prime}$ we see that $C^{\prime \prime}$ must be a bidegree $(a,0)$ curve through $Q_{2}, Q_{3},\ldots, Q_{a}, Q_{a+2}$, hence it must coincide with $F_{2}F_{3}\ldots F_{a}C_{2}$. Again, this contradicts the linear independence of $P_{1}$ and $P_{2}$.

\begin{figure}[h]
\begin{center} 
\scalebox{0.27}
{\includegraphics{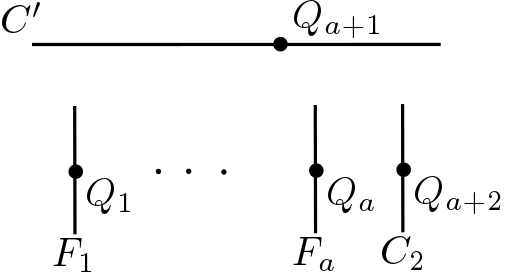}}
\caption{The configuration of $A$ when $m_{0,1}(A)=a+1$.}\label{fig10} 
\end{center}
\end{figure}

\textbf{Case 3:} Assume that $m_{0,1}(A)=a+2$. Say $C$ is a $(0,1)$-curve passing through all points of $A$. First, assume that $C$ is irreducible.   Since we assume that   $E_{\beta}^{+}(T)$ is not contained in a curve of bidegree $(0,1)$, there must exist a point $Q_{a+3}\in  E_{\beta}^{+}(T)$ which does not lie on $C$.  Since  no two points of $A$ lie on the same fiber of $\varphi$, $Q_{a+3}$ might be on the same fiber with at most one of the point of $A$, say $Q_{a+2}$.   Let $B=\{Q_{1}, Q_{2},\ldots, Q_{a+1}, Q_{a+3}\}$. By the choice of the set $B$,   no two points of $B$ lie on the same fiber of $\varphi$.  We claim that $m_{0,1}(B)=a+1$. Else, there exists a $(0,1)$-curve $\tilde{C}$ passing through all points of $B$. The intersection number of two $(0,1)$ curves on $\mathbb{F}_{a}$ is equal to $a$, but $C$ and $\tilde{C}$ share the $a+1$ points $Q_{1},\ldots, Q_{a+1}$. Since $C$ is irreducible, we must have $C=\tilde{C}$. This is a contradiction since $Q_{a+3}\notin C$. Therefore $m_{0,1}(B)=a+1$. We can now replace $A$ by $B$ and go back to Case 1 or Case 2. 

Now we assume that $C$ is reducible. Then we can write $C=C^{\prime} F_{1} F_{2}\ldots F_{a}$ where $C^{\prime}$ is the (unique, irreducible) curve of bidegree $(-a,1)$ and $F_{i}$'s are fibers of $\varphi$. By our initial assumption on the set $A$ that no two points of $A$ are on the same vertical fiber, $C'$ contains at least two points of $A$, say $Q_{a+1},Q_{a+2}\in C'$. Since we are assuming that  $E_{\beta}^{+}(T)$  is not contained in a $(0,1)$-curve, we can modify our choice of $A$ such that only $Q_{a+1}, Q_{a+2}$ are contained in $C^{\prime}$ ,  $Q_{1}, \ldots, Q_{a}$ are contained in $a$ distinct fibers $F_{1},\ldots, F_{a}$ and no two points of $A$ lie on the same fiber of $\varphi$.  Moreover, there exist a point $Q_{a+3}\in E^+_{\beta}(T)\setminus C$ (see Figure \ref{fig11}).  $Q_{a+3}$ might be on the same fiber with at most one of the point of $C'$, say $Q_{a+2}$.   Set $B=\{Q_{1}, Q_{2},\ldots, Q_{a+1}, Q_{a+3}\}$.   By the choice of the set $B$,   no two points of $B$ lie on the same fiber of $\varphi$.  We claim that $m_{0,1}(B)=a+1$. Else, there exists a $(0,1)$-curve $\tilde{C}$ passing through all points of $B$. The intersection number of $\tilde{C}$ and $C^{\prime}$ is $0$, but they share the point $Q_{a+1}$. Hence $C^{\prime}$ must divide $\tilde{C}$. But then $\tilde{C}$ must be divisible by all $F_{i}$ in order to contain $Q_{1},\ldots, Q_{a}$, so $\tilde{C}=C$. This is a contradiction since $Q_{a+3}$ is not on $C$. This finishes the proof of this case and the proof that $P_{1}$ and $P_{2}$ have no common factor, possibly after modifying the choice of $A$

%The remaining case is the following: Through any $a+2$ points in  $E_{\beta}^{+}(T)$  there passes a reducible $(0,1)$-curve. All such curves contain the unique irreducible $(-a,1)$-curve $C^{\prime}$ as a component along with $a$ fibers of $\varphi$. Since we are assuming that  $E_{\beta}^{+}(T)$  is not contained in a $(0,1)$-curve, we can modify our choice of $A$ such that only $Q_{a+1}, Q_{a+2}$ are contained in $C^{\prime}$ and $Q_{1}, \ldots, Q_{a}$ are contained in $a$ distinct fibers $F_{1},\ldots, F_{a}$. Again, set $C=C^{\prime}F_{1}\ldots F_{a}$ and choose $Q_{a+3}\in  E_{\beta}^{+}(T)$ not lying on $C$.  Set $B=\{Q_{1}, Q_{2},\ldots, Q_{a+1}, Q_{a+3}\}$.  Since $m_{0,1}(B)=a+2$ by assumption, there exists a $(0,1)$-curve $\tilde{C}$ passing through all points of $B$. The intersection number of $\tilde{C}$ and $C^{\prime}$ is $0$, but they share the point $Q_{a+1}$. Hence $C^{\prime}$ must divide $\tilde{C}$. But then $\tilde{C}$ must be divisible by all $F_{i}$ in order to contain $Q_{1},\ldots, Q_{a}$, so $\tilde{C}=C$. This is a contradiction since $Q_{a+3}$ is not on $C$. This finishes the proof of this case and the proof that $P_{1}$ and $P_{2}$ have no common factor, possibly after modifying the choice of $A$. 

\begin{figure}[h]
\begin{center} 
\scalebox{0.22}
{\includegraphics{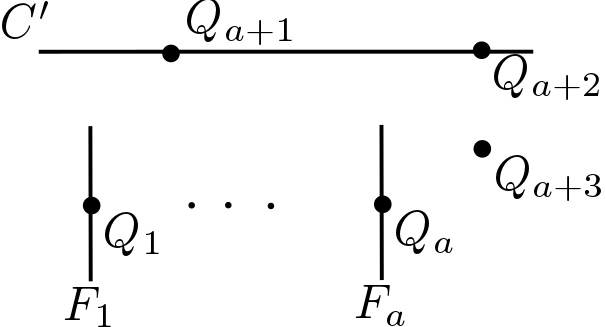}}
\caption{The configuration of $A$ when $m_{0,1}(A)=a+2$.}\label{fig11} 
\end{center}
\end{figure}

Now, By Theorem 3.1 in \cite{G}, $$\frac{1}{2}dd^{c}\log(|P_{1}|^2+|P_{2}|^2)$$ determines the $(1,1)$-current $S\in \mathcal T_{1,1}(\mathbb F_a).$

 Since $S$  is smooth except at finitely many logarithmic poles,  by a result of Demailly \cite{Dem},  the intersection $T\wedge S$ is a well-defined positive measure. We have 
 \[ \alpha+(a+1)\alpha^{\prime}=\int_{\mathbb{F}_{a}}T\wedge S  \geq \sum_{i=1}^{a+2} T\wedge S (Q_{i})   \geq \sum_{i=1}^{a+2} \nu(T, Q_{i})\nu(S, Q_{i}) >(a+2)\beta. \]
The second inequality follows from Demailly's comparison theorem for Lelong numbers  (Corollary 5.10 in \cite{Dem}). 
This gives us  $\beta< (\alpha + (a+1)\alpha^{\prime})/(a+2)$, which is a contradiction. This finishes the proof. 

\end{proof} 

\noindent \textbf{Acknowledgments.} We would like to thank the referee for his/her remarks and corrections which helped to improve the presentation of the paper. The first author is supported by CIMPA Research in Pairs Program and T\"UB\.ITAK 2219 Program. The second author is supported by T\"UB\.ITAK 3501 Proj. No 120F084  and T\"UB\.ITAK 2518 Proj. No. 119N642.  \\

\noindent \textbf{Conflict of interest.} The authors declare that they have no conflict of interest.

\end{document}